\documentclass[a4paper, 12pt]{article}

\usepackage[koi8-r,cp1251]{inputenc}
\usepackage{amsfonts,amssymb,amsmath, hyperref}
\usepackage[final]{epsfig}
\usepackage{graphicx}

\begin{document}

\begin{center}
\textbf{ THE STRUCTURE  OF NORMAL HAUSDORFF OPERATORS ON LEBESGUE SPACES}
\end{center}

\begin{center}
\textsc{A. R. Mirotin}
\end{center}

\begin{center}
 amirotin@yandex.ru
\end{center}

\

\textsc{Abstract.} {\small We consider a generalization of Hausdorff operator and  introduce the notion of the symbol of such an operator. Using this notion we describe  the   structure and investigate important properties  (such as invertibility, spectrum, norm, and compactness) of normal generalized  Hausdorff operators  on Lebesgue spaces over $\mathbb{R}^n.$ The examples of Ces\`{a}ro operators are considered.}

\

\textsc{Introduction.} Hausdorff operators on Lebesgue spaces   originated from some classical summation methods. They were introduced
by Hardy \cite[Chapter XI]{H} on the unit interval and by
 C. Georgakis and independently by
E. Liflyand and F. Moricz on the whole real line  \cite{G}, \cite{LM}.
  Its natural multidimensional extension  is the operator
$$
(\mathcal{H}f)(x) =\int_{\mathbb{R}^m} \Phi(u)f(A(u)x)du,\  x\in \mathbb{R}^n
$$
where  $A(u)$ is a family of $n\times n$ matrices satisfying $\mathrm{det}A(u) \ne 0$ almost everywhere in the support of  $\Phi.$ It was introduced by Brown and Moricz \cite{BM} on the Lebesgue space and  by Lerner and Liflyand \cite{LL} on the real
Hardy space.
 This class of operators has  attracted considerable attention in recent decades. It includes  some important examples, such as the Hardy operator, the adjoint
Hardy operator,  the Ces\`{a}ro operator and their multidimensional analogs.  The survey articles   \cite{Ls}, \cite{CFW}   contain
main results on Hausdorff operators and bibliography up to 2014. As far as the author is aware all known results refer  to the boundedness of general Hausdorff operators in various settings only. But this question, being solved positively, naturally entails such questions as whether the operator is invertible and what is its inverse, what is the spectrum of the operator, and a number of others. To deepened the investigations it is natural to begin with Hausdorff operators on the Hilbert space. This work  is mainly devoted to this case.
We shall consider a generalization of Hausdorff operator and  introduce the notion of the symbol of this operator. This notion is crucial for our considerations. Our aim  is to describe  the structure and to investigate  important properties  (such as invertibility, spectrum, norm, and compactness) of generalized normal Hausdorff operators  on Lebesgue spaces over $\mathbb{R}^n.$  Using the Mellin transform we prove under some extra condition that a normal Hausdorff operator in $L^2(\mathbb{R}^n)$ is unitary equivalent to the operator   of coordinate-wise multiplication by  the symbol in the space $L^2(\mathbb{R}^n,\mathbb{C}^{2^n}).$
As  examples,  discrete   Hausdorff operators and   Ces\`{a}ro operators are considered. 
 The study of the one-dimensional  Ces\`{a}ro operator on $L^2$ was pioneered by Brown, Halmos, and  Shields \cite{BHS}, the $L^p$ case was considered in \cite{BM}.

\

\textbf{Definition 1 } \cite{arx}. Let  $(\Omega,\mu)$ be some $\sigma$-compact locally compact topological space endowed with positive regular  Borel measure $\mu,$ $\Phi$  a locally integrable function on $\Omega,$  and $(A(u))_{u\in \Omega}$  a $\mu$-measurable family of $n\times n$ matrices defined  almost everywhere in  the support of  $\Phi$ and satisfying $\mathrm{det}A(u) \ne 0.$  We define the \textit{Hausdorff  operator} with the kernel $\Phi$  by
$$
(\mathcal{H}f)(x) = (\mathcal{H}_{\Phi, A}f)(x) =\int_\Omega \Phi(u)f(A(u)x)d\mu(u)
$$
($x\in\mathbb{ R}^n$ is a column vector).

\

As was mentioned by Hardy in the case $n=1,$   $\Omega=[0,1]$ \cite[Theorem 217]{H} the Hausdorff  operator possesses some regularity property. The multidimensional version of his result looks as follows.

\textbf{Proposition 1.} \textit{Under the conditions of definition 1 let the matrices $A(u)$ be positive definite for  almost all $u\in \Omega.$  In order that the transformation $\mathcal{H}_{\Phi, A}$  should be
regular, i.e. that $f\in C(\mathbb{R}^n),$ $f(x) \to l$ when $x\to \infty$   should imply $\mathcal{H}_{\Phi, A}f(x) \to l$, it is necessary and
sufficient that $\int_\Omega \Phi(u)d\mu(u)=1.$
}

Proof. If $f(x)=1$ then $\mathcal{H}f(x) =\int_\Omega \Phi(u)d\mu(u).$ Thus, $\int_\Omega \Phi(u)d\mu(u)=1$ is a necessary condition.

To prove the sufficiency, first note that if $A(u)$ is positive definite then  $x\to \infty$ implies $A(u)x\to \infty$ (this follows from the representation $A(u)=C^{-1}A'C$ where $C$ is an orthogonal matrix and $A'$ is a diagonal one with positive eigenvalues). But if, in addition,  $f\in C(\mathbb{R}^n)$ and $f(x) \to l$ then $f$ is bounded and therefore  $\mathcal{H}f(x) \to l$ ($x\to \infty$) by the Lebesgue Theorem.

Proposition 1 motivates our work with  positive definite matrices in the sequel.
 
\

\textsc{1. The structure of normal Hausdorff  operators in $L^2.$}

\

We need the following lemmas to prove our main result.

\

\textbf{Lemma 1} \cite{arx} (cf. \cite[(11.18.4)]{H}, \cite{BM}). \textit{Let  $(\det A(u))^{-1/p}\Phi(u)\in L^1(\Omega).$ Then the operator $\mathcal{H}_{\Phi, A}$ is bounded in $L^p(\mathbb{R}^n)$ ($1\leq p<\infty$) and }
$$
\|\mathcal{H}_{\Phi, A}\|\leq \int_\Omega |\Phi(u)||\det A(u)|^{-1/p}d\mu(u).
$$

This estimate is sharp (see corollary 2 below).

\

In the next corollary we as usual consider $L^p(\mathbb{R}_+^n)$ as a subspace of $L^p(\mathbb{R}^n).$

\textbf{Corollary 1.} \textit{If, in addition, every $A(u)$ maps $\mathbb{R}_+^n$ into itself the operator $\mathcal{H}_{\Phi, A}$ is bounded in $L^p(\mathbb{R}_+^n)$ ($1\leq p<\infty$).}

\

\textbf{Lemma 2} (cf. \cite{BM}). \textit{The adjoint for the Hausdorff operator in $L^p(\mathbb{R}^n)$ has the form}
$$
(\mathcal{H}^*f)(x) = (\mathcal{H}_{\Phi, A}^*f)(x) =\int_\Omega \overline{\Phi(v)}|\det A(v)|^{-1}f(A(v)^{-1}x)d\mu(v).
$$
Thus, the adjoint for a Hausdorff operator  is also  Hausdorff.

\

\textbf{Lemma 3}.  \textit{Under the conditions of definition 1 the Hausdorff operator $\mathcal{H}_{\Phi, A}$ is normal in $L^2(\mathbb{R}^n)$ if   $n\times n$-matrices $A(u)$ (defined $\mu$-a.~e. on the support of $\Phi$) form a  pairwise commuting family.
}

Proof. Indeed, this follows from the equalities
$$
(\mathcal{H}\mathcal{H}^*f)(x)=\int_\Omega \int_\Omega\Phi(u)\overline{\Phi(v)}|\det A(v)|^{-1}f(A(u)A(v)^{-1}x)d\mu(v)d\mu(u),
$$
$$
(\mathcal{H}^*\mathcal{H}f)(x)=\int_\Omega \int_\Omega\Phi(u)\overline{\Phi(v)}|\det A(v)|^{-1}f(A(v)^{-1}A(u)x)d\mu(u)d\mu(v)
$$
and the Fubini Theorem.

\

In the sequel for $x\in\mathbb{ R}^n, \alpha\in \mathbb{C}^n,$ and $r\in \mathbb{R}$ we put  $dx=$ $dx_1\dots dx_n,$ $|x|^{r+\alpha}:=$ $|x_1|^{r+\alpha_1}\dots |x_n|^{r+\alpha_n}.$

\
The following notion is crucial for our investigations.

\textbf{Definition 2.} Let the conditions of definition 1 are fulfilled. Let $a(u):=$ $(a_1(u),\dots,a_n(u))$ be the family of eigenvalues  (with their multiplicities) of the matrix $A(u)$  and  $|\det A(u)|^{-1/p}\Phi(u)\in L^1(\Omega).$  Then we call the function $(s\in\mathbb{ R}^n)$
   $$
   \varphi(s):=\int_{\Omega}\Phi(u)(\det A(u))^{-1/p}a(u)^{-is}d\mu(u)
   $$
\textit{the symbol of the Hausdorff operator} $\mathcal{H}_{\Phi, A}$ in $L^p(\mathbb{R}^n)$ ($1\leq p<\infty$).

In the previous definition we assume that for every $j=1,\dots, n$  a $\mu$-measurable branch of $\log a_j(u)$ in $\Omega$ exists and is fixed, and  $a(u)^{z}:=$ $\prod_{j=1}^n a_j(u)^{z_j}$ where $a_j(u)^{z_j}:=$ $\exp(z_j\log a_j(u))$ ($z\in \mathbb{C}^n$). (Otherwise, $ \varphi$ should be considered as a multi-valued function.)

It is clear that $\varphi$ is bounded and continuous on $\mathbb{R}^n.$

\

Note that if $A(u)$ is a pairwise commuting family of real selfadjoint $n\times n$-matrices, there is such an orthogonal $n\times n$-matrix $C$ and a family of real diagonal  matrices $A'(u)=\mathrm{diag}[a_1(u),\dots,a_n(u)],$ that $A'(u)=C^{-1}A(u)C$ for  almost all  $u\in \Omega.$
 If each $A(u)$ is positive definite the corresponding Hausdorff operator $\mathcal{H}_{\Phi, A'}$  in $L^p(\mathbb{R}^n)$ is bounded provided   $(\det A(u))^{-1/p}\Phi(u)=a(u)^{-1/p}\Phi(u)\in L^1(\Omega).$

\

\textbf{Example 1} (cf. \cite{BM}). Consider the \textit{$n$-dimensional  Ces\`{a}ro operator}:
$$
(\mathcal{C}_nf)(x_1,\dots, x_n)=\frac{1}{x_1 \dots x_n}\int_0^{x_1}\!\!\dots\int_0^{x_n}f(t_1,\dots, t_n)dt_1\dots dt_n.
$$
This is a bounded  Hausdorff operator in $L^p(\mathbb{R}^n)$ (and in $L^p(\mathbb{R}_+^n)$) for $1<p<\infty$ where $\Omega=[0,1]^n$ is endowed with the Lebesgue measure, $\Phi=1$, and $A(u)=\mathrm{diag}[u_1,\dots, u_n].$ Its symbol is ($s=(s_j)\in\mathbb{ R}^n$)
$$
\varphi(s)=\int_0^1\!\!\dots\int_0^1 \prod\limits_{j=1}^n u_j^{-\frac{1}{2}-is_j}du_1\dots du_n=\prod_{j=1}^n(1/2-is_j)^{-1}.
$$

\

\textbf{Example 2}. (Discrete   Hausdorff operators.) Let $\Omega=\mathbb{Z}_+,$ and $\mu$ be a counting measure. Then the definition 1 takes the form ($f\in L^p(\mathbb{R}^n)$)
$$
(\mathcal{H}_{\Phi, A}f)(x) =\sum_{k=0}^\infty \Phi(k)f(A(k)x).
$$
Assume that $\sum_{k=0}^\infty |\Phi(k)||\det A(k)|^{-1/p}<\infty.$ Then $\mathcal{H}_{\Phi, A}$ is bounded on $L^p(\mathbb{R}^n)$ and  if $a_j(k)>0$ for all $j$ and $k$  the symbol of $\mathcal{H}_{\Phi, A}$ is
$$
\varphi(s)=\sum_{k=0}^\infty \Phi(k)(\det A(k))^{-1/p}a(k)^{-is}
$$
where $s=(s_j)\in\mathbb{ R}^n$ and the principal values of the exponents  are considered.
Since this series converges  on $\mathbb{R}^n$ absolutely (and therefore uniformly), $\varphi$ is uniformly almost periodic.

\

 The next theorem   describes under some additional assumptions the structure of  normal Hausdorff operator in $L^2(\mathbb{R}^n).$

\textbf{Theorem 1.} \textit{Let $A(u)$ be a pairwise commuting family of real  positive definite $n\times n$-matrices  (defined $\mu$-a.~e. on the support of $\Phi$), and $(\det 
A(u))^{-1/2}\Phi(u)\in L^1(\Omega).$ Then the  Hausdorff operator $\mathcal{H}_{\Phi, A}$ in $L^2(\mathbb{R}^n)$  with symbol $\varphi$ is normal and   unitary equivalent to the operator $M_\varphi$  of coordinate-wise multiplication by  $\varphi$ in the space $L^2(\mathbb{R}^n,\mathbb{C}^{2^n})$ of $\mathbb{C}^{2^n}$-valued functions. More precisely,  $\mathcal{H}_{\Phi, A}=\mathcal{V}^{-1}M_\varphi \mathcal{V}$ where  $\mathcal{V}$ is a unitary operator from $L^2(\mathbb{R}^n)$ to $L^2(\mathbb{R}^n,\mathbb{C}^{2^n})$ which does not depend on $\Phi.$}

\

Proof. Since matrices $A(u)$  pairwise commutes and are positive definite, there are an   orthogonal $n\times n$-matrix $C$ and a family of diagonal  positive definite matrices $A'(u)$ such that $A'(u)=C^{-1}A(u)C$ for  all $u\in \Omega.$ If $A'(u)=\mathrm{diag}[a_1(u),\dots,a_n(u)],$  then all functions $a_j(u)$ are $\mu$-measurable and positive and $\det A(u)=$ $a_1(u)\dots a_n(u).$ Consider the unitary operator $\hat{C}f(x):=f(Cx)$
in $L^2(\mathbb{R}^n).$ It is easy to verify that
$$
\mathcal{H}_{\Phi, A}=\hat{C}^{-1}\mathcal{H}_{\Phi, A'}\hat{C}.
$$
Thus, the  operator $\mathcal{H}_{\Phi, A}$  is unitary equivalent to $\mathcal{H}_{\Phi, A'}.$

 Let $U_1, \dots, U_{2^n}$ be open hyperoctants in $\mathbb{R}^n.$ Then $L^2(\mathbb{R}^n)=\bigoplus_{j=1}^{2^n}L^2(U_j),$ the orthogonal sum of subspaces, and each $L^2(U_j)$ is $\mathcal{H}_{\Phi, A'}$-invariant. For every $j$ consider the operator $\mathcal{H}_j'$ on $L^2(\mathbb{R}^n)$ which on $L^2(U_j)$  coincides with the restriction of $\mathcal{H}_{\Phi, A'}$ to this space and equals to zero on the orthogonal compliment of $L^2(U_j).$  Then
$$
\mathcal{H}_{\Phi, A'}=\bigoplus_{j=1}^{2^n}\mathcal{H}_j',
$$
the orthogonal sum of operators.

 Consider the modified  $n$-dimensional  Mellin transform for the $n$-hyperoctant $U_j$  in the form ($j=1,\dots,2^n;\ s\in \mathbb{R}^n$)
$$
(\mathcal{M}_jf)(s):=\frac{1}{(2\pi)^{n/2}}\int_{U_j}|x|^{-\frac{1}{2}+is}f(x)dx.
$$
Then $\mathcal{M}_j$ is  a unitary operator from $L^2(U_j)$ to $L^2(\mathbb{R}^n).$ This can be easily obtained from the Plancherel theorem  for the $n$-dimensional Fourier transform
by using an exponential change of variables  (see \cite{BPT}). Moreover, if we assume that   $|y|^{-1/2}f(y)\in L^1(U)$ then making use of Fubini’s theorem, and integrating
by substitution $x=A'(u)^{-1}y$ yield the following ($s\in \mathbb{R}^n$):
$$
(\mathcal{M}_j \mathcal{H}_j'f)(s)=\frac{1}{(2\pi)^{n/2}}\int_{U_j}|x|^{-\frac{1}{2}+is}dx\int_{\Omega}\Phi(u)f(A'(u)x)d\mu(u)=
$$
$$
\frac{1}{(2\pi)^{n/2}}\int_{\Omega}\Phi(u)d\mu(u)\int_{U_j}|x|^{-\frac{1}{2}+is}f(A'(u)x)dx=
$$
$$
\int_{\Omega}\Phi(u)a(u)^{-\frac{1}{p}-is}d\mu(u)\frac{1}{(2\pi)^{n/2}}\int_{U}|y|^{-\frac{1}{2}+is}f(y)dy=
\varphi(s)(\mathcal{M}_jf)(s).
$$
By continuity we get for all $f\in L^2(U_j)$ that
$$
\mathcal{M}_j\mathcal{H}_j'f=\varphi \mathcal{M}_jf.
$$
So, $\mathcal{H}_j'$ is unitary equivalent to the operator $M'_\varphi$ in $L^2(\mathbb{R}^n)$ of multiplication by $\varphi$ and
$$
\mathcal{H}_{\Phi, A'}=\bigoplus_{j=1}^{2^n}\mathcal{M}_j^{-1}M'_\varphi\mathcal{M}_j.
$$
If we define the operator $\mathcal{V}'$ from $L^2(\mathbb{R}^n)$ to the space $L^2(\mathbb{R}^n)\oplus\dots\oplus L^2(\mathbb{R}^n)$ ($2^n$ direct summands) which is isomorphic to the Hilbert space $L^2(\mathbb{R}^n,\mathbb{C}^{2^n})$  by the equalities $\mathcal{V}'|L^2(U_j):=\mathcal{M}_j$ ($j=1,\dots,2^n$) then $\mathcal{V}'$ is unitary and it is easy to verify that
$$
\mathcal{V}'^{-1}M_\varphi \mathcal{V}'=\bigoplus_{j=1}^{2^n}\mathcal{M}_j^{-1}M'_\varphi\mathcal{M}_j.
$$
Thus, $\mathcal{H}_{\Phi, A'}=\mathcal{V}'^{-1}M'_\varphi \mathcal{V}'$ and the proof is complete.

\

The following corollary shows that properties of a Hausdorff operator are closely related to the properties of its symbol.

\textbf{Corollary 2.}
(i)  \textit{The   operator $\mathcal{H}_{\Phi, A}$ is invertible if and only if $\inf |\varphi|>0;$
in this case the inverse $\mathcal{H}^{-1}_{\Phi, A}$
is  unitary equivalent to the   operator  $M_{1/\varphi}$ in $L^2(\mathbb{R}^n,\mathbb{C}^{2^n})$;}

(ii) \textit{Let $\mathcal{H}_{\Phi, A}$ and $\mathcal{H}_{\Psi, B}$ be two  Hausdorff operators over the same measure space $(\Omega, \mu)$ such that $(A(u), B(v))$ is a pairwise commuting family of real  positive definite $n\times n$-matrices  ($u$ and $v$ run over the support of $\Phi$ and $\Psi$ respectively), and $(\det A(u))^{-1/2}\Phi(u)$, $(\det B(v))^{-1/2}\Psi(v) \in L^1(\Omega).$ Then the product  $\mathcal{H}_{\Phi, A}\mathcal{H}_{\Psi, B}$ is  unitary equivalent to the operator $M_{\varphi\psi}$ in $L^2(\mathbb{R}^n,\mathbb{C}^{2^n})$ ($\psi$ denotes the symbol of $\mathcal{H}_{\Psi, B}$).}

(iii)  \textit{The   spectrum $\sigma(\mathcal{H}_{\Phi, A}),$  the  point spectrum $\sigma_p(\mathcal{H}_{\Phi, A}),$ and  the continuous  spectrum of $\mathcal{H}_{\Phi, A}$  equal to the spectrum (i.~e. to the closure of the range of  the   symbol $\varphi$), to the  point spectrum, and to the continuous  spectrum  of $M'_{\varphi}$ in $L^2(\mathbb{R}^n)$ respectively, the residual spectrum $\sigma_r(\mathcal{H}_{\Phi, A})$ is empty};

(iv) \textit{The symbol  of the adjoint $\mathcal{H}_{\Phi, A}^*$ is $\overline{\varphi}$.}

(v) $\|\mathcal{H}_{\Phi, A}\|=\sup |\varphi|.$ \textit{If, in addition,  $\Phi(u)\geq 0$ for $\mu$-almost all $u$ then}
$$
\|\mathcal{H}_{\Phi, A}\|=\int_\Omega \Phi(u)(\det A(u))^{-1/2}d\mu(u).
$$

Proof. (i) Evidently  the   operator $\mathcal{H}_{\Phi, A}$ is invertible if and only if $M_\varphi$ is invertible, i.~e. if $\inf |\varphi|>0,$ and in this case  $\mathcal{H}^{-1}_{\Phi, A}$
is  unitary equivalent to  $M_{1/\varphi}.$

To prove (ii) first note that the orthogonal matrix $C$ exists such that both  $A'(u)=C^{-1}A(u)C$ and $B'(v):=C^{-1}B(v)C$ are diagonal. Then the proof of theorem 1 shows that $\mathcal{V}\mathcal{H}_{\Phi, A}\mathcal{V}^{-1}=M_\varphi ,$ $ \mathcal{V}\mathcal{H}_{\Psi, B}\mathcal{V}^{-1}=M_\psi$ for some unitary operator $\mathcal{V}$ from $L^2(\mathbb{R}^n)$ to $L^2(\mathbb{R}^n,\mathbb{C}^{2^n})$ and (ii) follows.

(iii) For $\lambda\in \mathbb{C}$ we have
$$
M_\varphi-\lambda=\bigoplus_{1}^{2^n}(M'_\varphi-\lambda):=(M'_\varphi-\lambda)\oplus\dots\oplus (M'_\varphi-\lambda).
$$
 Since a finite
 orthogonal sum of operators is invertible  if and only if each summand is invertible   (see, e.g.,  \cite[p. 439]{Kub}),  $\sigma(\mathcal{H})=\sigma(M_\varphi)=\sigma(M'_\varphi).$

Moreover,
 $$
 \sigma_p(\mathcal{H})=\sigma_p(M_\varphi)=\sigma_p\left(\bigoplus_{1}^{2^n} M'_\varphi\right)=
 \bigcup_{1}^{2^n}\sigma_p(M'_\varphi)=\sigma_p(M'_\varphi).
 $$

  Next, by theorem 1 $\sigma_r(\mathcal{H})=\sigma_r(M_\varphi).$  Let $\lambda\notin  \sigma_p(M_\varphi).$  Since $\sigma_r(M'_\varphi)=\emptyset,$ we get (in the following $\mathrm{cl}(S)$ denotes the closure of the subset $S\subset L^2(\mathbb{R}^n)$ and $\mathrm{Im}(T)$ denotes the image of the operator $T$)
  $$
 \mathrm{cl}(\mathrm{Im}(M_\varphi-\lambda)) =\bigoplus_{1}^{2^n} \mathrm{cl}(\mathrm{Im}(M'_\varphi-\lambda))=\bigoplus_{1}^{2^n}L^2(\mathbb{R}^n)=L^2(\mathbb{R}^n,\mathbb{ C}^{2^n}).
 $$
So, $\sigma_r(M_\varphi)=\emptyset,$   and (iii) follows.

(iv)  Since the adjoint of $\mathcal{H}_{\Phi, A}$ is also of Hausdorff type (with $A(u)^{-1}$ instead of $A(u)$), the this statement follows  from the definition 2 and lemma 2.

Finally, the equality $\|\mathcal{H}_{\Phi, A}\|=\sup |\varphi|$
 follows  from (iii) and the normality of $\mathcal{H}_{\Phi, A}$. If, in addition,  $\Phi(u)\geq 0$ for $\mu$-almost all $u$ we have $\sup |\varphi|=\varphi(0)$ which   completes the proof.

\

\textbf{Corollary 3.} \textit{Under the assumptions of theorem 1 the  Hausdorff operator $\mathcal{H}_{\Phi, A}$ is self-adjoint (positive, unitary) in $L^2(\mathbb{R}^n)$ if and only if its symbol $\varphi$ is real-valued (respectively, nonnegative, $|\varphi|=1$).} 

Proof. Indeed, the spectral decomposition
of a normal operator implies that  the operator $\mathcal{H}_{\Phi, A}$ is self-adjoint (positive, unitary) if and only if $\sigma(\mathcal{H}_{\Phi, A})=\mathrm{cl}(\varphi(\mathbb{R}^n))$ is contained in $\mathbb{R}$ (respectively, in $\mathbb{R}_+,$  $\mathbb{T}$).

\

\textbf{Corollary 4.} \textit{Under the assumptions and notation of theorem 1 (ii)  the  operator $\mathcal{H}_{\Psi, B}$ is the  inverse of $\mathcal{H}_{\Phi, A}$ if and only if $\varphi\psi=1$.}

\

\textbf{Corollary 5.} \textit{Under the assumptions of theorem 1 the  Hausdorff operator $\mathcal{H}_{\Phi, A}$ in $L^2(\mathbb{R}_+^n)$ is unitary equivalent to the operator of multiplication by  $\varphi$ in $L^2\mathbb{(R}_+^n).$}

Proof. This was established in the process of proving theorem 1.

\

\textbf{Example 3}. Consider the  Ces\`{a}ro operator  in $L^2(\mathbb{R}^n)$ (see example 1).
It is normal and its symbol is $\varphi(s)=\prod_{j=1}^n(1/2-is_j)^{-1}.$

The statement (iv) of corollary 2 implies $\|\mathcal{C}_n\|=2^n$ (and by the corollary 5 the  same is true in $L^2(\mathbb{R}_+^n)$).

 According to the statement (iii) of corollary 2 $\sigma_r(\mathcal{C}_n)=\emptyset,$ $\sigma_p(\mathcal{C}_n)=\emptyset,$ and $\sigma(\mathcal{C}_n)$ equals to the closure of the range of $\varphi.$ Consider the circle $S:=\{1/(1/2-it): t\in \mathbb{R}\} =\{z\in \mathbb{C}: |z-1|=1\}.$ Then  is $ \sigma(\mathcal{C}_n)=\{z_1\dots z_n: z_j\in S,j=1,\dots, n\}.$
 Let $z_j\in S,$  $z_j= 1+e^{i\theta_j}$ ($\theta_j\in [-\pi,\pi]$). Then $|z_j|=2\cos(\theta_j/2),$ and $\arg (z_j)=\theta_j/2$ $(j=1,\dots,n).$ Let $z=z_1\dots z_n,$ $r:=|z|,$ $\theta:=\arg(z).$ Then $\theta=(\theta_1+\dots+\theta_n)/2,$ and $r=2^n\prod_{j=1}^n\cos(\theta_j/2).$
Using the identity
 $$
 \prod_{j=1}^n\cos(\alpha_j) =\frac{1}{2^{n-1}}\sum\limits_{\varepsilon\in\{-1,1\}^n}
 \cos(\varepsilon\cdot\alpha)=\frac{1}{2^{n-1}}\sum\limits_{\varepsilon'\in\{-1,1\}^{n-1}}
 \cos(\alpha_1+\varepsilon'\cdot\alpha'),
 $$
 where $\alpha= (\alpha_j)_{j=1}^n=$ $(\theta_j/2)_{j=1}^n,$ $\alpha'=(\theta_j/2)_{j=2}^n,$ and $\cdot$ denotes the dot product we conclude that
 $$
\sigma(\mathcal{C}_n)=\{re^{i\theta}: r\leq 2(2^{n-1}-1+\cos \theta), \theta\in [-\pi,\pi]\}.
$$
 By the corollary 5 the  same is true for  the  Ces\`{a}ro operator  in $L^2(\mathbb{R}_+^n).$

  In particular, for  the bivariate Ces\`{a}ro operator we have
$$
\sigma(\mathcal{C}_2)=\{re^{i\theta}: r\leq 2(1+\cos \theta), \theta\in [-\pi,\pi]\},
$$
the region in $\mathbb{C}$ bounded by  cardioid (cf. \cite[p. 95]{R}).

Results of the example 3 are consistent with classical results from  \cite{BHS}  where the case $n=1$ and the space $L^2(\mathbb{R}_+)$ were considered (see p. 137 therein; in \cite{BHS} the Ces\`{a}ro operator  in $L^2(\mathbb{R}_+)$ was denoted by $\mathcal{C}_\infty$).

\

\textbf{Example 4}. Consider the $(C, k)$ mean of a function $f\in L^2(\mathbb{R})$ (see, e.g., \cite[p. 276]{H})
$$
(C, k)f(x)=\frac{k}{x^k}\int_0^x(x-t)^{k-1}f(t)dt\ (k>0).
$$
This is an operator of Hausdorff type where $n=1,$ $\Omega=[0,1]$ endowed with the Lebesgue measure, $\Phi(u)=k(1-u)^{k-1},$ and $A(u)=u.$ Its symbol is
$$
\varphi(s)=\int_0^1 k(1-u)^{k-1}u^{-1/2-is}du=\frac{\Gamma(k+1)\Gamma(1/2-is)}{\Gamma(k+1/2-is)}.
$$
Using \cite[Section 12.13, Example 1]{WW} we can rewrite this as follows
$$
\varphi(s)=\prod_{l=1}^\infty\frac{l(k+l-1/2-is)}{(k+l)(l-1/2-is)}.
$$
So, the spectrum  $\sigma((C,k))$ is the closure of the curve with complex equation 
$$
z=\prod_{l=1}^\infty\frac{l(k+l-1/2-is)}{(k+l)(l-1/2-is)},\ s\in \mathbb{R}.
$$
Moreover, the last representation of $\varphi$ implies also that $\max_{\mathbb{R}}|\varphi|=\varphi(0).$
It  follows in view of corollary 2 that
$$
\|(C,k)\|=\sqrt{\pi}\frac{\Gamma(k+1)}{\Gamma(k+1/2)}
$$ 
(since $\Phi\geq 0,$ this follows
also from the lust equality in  corollary 2 (v)).

If, in addition, $k\in \mathbb{N},$ we have $\|(C,k)\|=k!2^k/1\cdot 3\dots (2k-1)$ and
$$
\varphi(s)=k!\prod_{j=0}^{k-1}(j+1/2-is)^{-1}.
$$

\

\textbf{Example 5}. Let $\mathcal{H}_{\Phi, A}$ and $\varphi$  be as in the example 2. If $p=2$ corollary 2 implies that $\sigma(\mathcal{H}_{\Phi, A})=\mathrm{cl}(\varphi(\mathbb{R}^n)).$ Assume, in addition, that $A(k)=A^k$ where $A\ne E$ be an $n\times n$-matrix with positive eigenvalues $\lambda_1,\dots,\lambda_n$ (taking into account their multiplicities). Consider the function ($z\in\overline{\mathbb{D}}:=\{|z|\leq 1$\})
$$
F(z):=\sum\limits_{k=0}^\infty \Phi(k)(\det A)^{-k/2}z^k.
$$
Then $F$ belongs to the commutative Banach algebra $\mathcal{A}^+$ of functions on $\overline{\mathbb{D}}$ with absolutely convergent Taylor series and $\varphi(s)=F(\exp(-i\sum_js_j\log \lambda_j)).$ Now theorem 1 implies that $\sigma(\mathcal{H}_{\Phi, A})=F(\mathbb{T})$ where $\mathbb{T}=\{|z|=1\}.$ It follows also that $\|\mathcal{H}_{\Phi, A}\|=\sup_{\mathbb{T}}|F|,$ and  $\mathcal{H}_{\Phi, A}$ is invertible if and only if $\inf_{\mathbb{T}}|F|>0.$ Since the last condition is equivalent to  $\inf_{\overline{\mathbb{D}}}|F|>0,$ under this  condition  the function $G:=1/F$ belongs to the algebra $\mathcal{A}^+,$ as well (it follows from the Gelfand theory). Let $G(z)=\sum_{k=0}^\infty b(k)z^k$ be the Taylor expansion  for $G$ ($z\in \overline{\mathbb{D}}$). Then the inverse of  $\mathcal{H}_{\Phi, A}$ is the Hausdorff operator
$$
\mathcal{H}_{\Phi, A}^{-1}f(x):=\sum\limits_{k=0}^\infty b(k)(\det A)^{k/2}f(A^kx).
$$
Indeed, in this  case all the conditions of corollary 2 (ii) for the pair  $\mathcal{H}_{\Phi, A},$ $\mathcal{H}_{\Phi, A}^{-1}$ are fulfilled and the symbol of the
last operator is
$$
\psi(s)=\sum\limits_{k=0}^\infty b(k)a(k)^{-is} =
$$
$$
G(\exp(-i\sum_js_j\log \lambda_j))=\frac{1}{F(\exp(-i\sum_js_j\log \lambda_j))}=\frac{1}{\varphi(s)}.
$$
Thus, the result follows from corollary 4.

\

\textsc{2. The noncompactness of Hausdorff  operators in $L^p.$} The problem of compactness of (nontrivial) Hausdorff operators  was posed by Liflyand in \cite{L} (see also \cite{Ls}). Below we shall shaw that
such an operator  is noncompact in $L^p(\mathbb{R}^n)$ ($1\leq p<\infty$) if $(A(u))_{u\in \Omega}$ is a pairwise commuting family of   positive definite matrices. The case $p=2$ was considered in \cite{arx}.

\textbf{Theorem 2.}
\textit{Let $A(u)$ be a pairwise commuting family of  positive definite $n\times n$-matrices ($u$ runs over the support of $\Phi$) and $(\det A(u))^{-1/p}\Phi(u)\in L^1(\Omega).$ Then the  Hausdorff operator $\mathcal{H}_{\Phi, A}$ is  noncompact in $L^p(\mathbb{R}^n)$ ($1\leq p<\infty$) provided it is non-zero.}

\

Proof. Assume the contrary. As in the proof of theorem 1 consider the (bounded and invertible) operator $\hat{C}f(x):=f(Cx)$
in $L^p(\mathbb{R}^n).$ Since $\hat{C}\mathcal{H}_{\Phi, A}\hat{C}^{-1}=\mathcal{H}_{\Phi, A'},$ the operator $\mathcal{H}:=\mathcal{H}_{\Phi, A'}$ is  non-zero and compact, too. Note that  every open $n$-hyperoctant in $\mathbb{R}^n$ is $\mathcal{H}$-invariant.
There is an open $n$-hyperoctant  $U$ such that the restriction $\mathcal{K}:=\mathcal{H}|L^p(U)$ is non-zero. Then  $\mathcal{K}$ is  compact as well.

We shall use the modified  $n$-dimensional  Mellin transform for the $n$-hyperoctant $U$  in the form ($s\in \mathbb{R}^n$)
$$
(\mathcal{M}f)(s):=\frac{1}{(2\pi)^{n/2}}\int_{U}|x|^{-\frac{1}{q}+is}f(x)dx
$$
to get a contradiction. The map $\mathcal{M}$ is a bounded operator from $L^p(U)$ to $L^q(\mathbb{R}^n)$ if $1\leq p\leq 2$  ($1/p+1/q=1$) and a unitary operator from $L^2(U)$ to $L^2(\mathbb{R}^n).$ This can be easily obtained from the Hausdorff–Young inequality (and from the Plancherel theorem if $p=2$) for the $n$-dimensional Fourier transform
by using an exponential change of variables
 (see \cite{BPT}). Let $f\in L^p(U).$  First assume that   $|y|^{-1/q}f(y)\in L^1(U).$ Then as in the proof of theorem 1 making use of Fubini’s theorem, and integrating
by substitution $x=A(u)^{-1}y$ yield the following ($s\in \mathbb{R}^n$):
$$
(\mathcal{MK}f)(s)=\frac{1}{(2\pi)^{n/2}}\int_{U}|x|^{-\frac{1}{q}+is}dx\int_{\Omega}\Phi(u)f(A(u)x)d\mu(u)=
$$
$$
\frac{1}{(2\pi)^{n/2}}\int_{\Omega}\Phi(u)d\mu(u)\int_{U}|x|^{-\frac{1}{q}+is}f(A(u)x)dx=
$$
$$
\int_{\Omega}\Phi(u)a(u)^{-\frac{1}{p}-is}d\mu(u)\frac{1}{(2\pi)^{n/2}}\int_{U}|y|^{-\frac{1}{q}+is}f(y)dy=
\varphi(s)(\mathcal{M}f)(s),
$$
where the symbol $\varphi$ is continuous and bounded on   $\mathbb{R}^n.$
Thus, if $|y|^{-1/q}f(y)\in L^1(U)$ then
$$
\mathcal{MK}f=\varphi \mathcal{M}f. \eqno(1)
$$
 By continuity the last equality  is valid for all $f\in L^p(U).$

 If $p=2$ it follows that $\mathcal{K}$ is unitary equivalent to the operator $M_\varphi$  of multiplication by $\varphi$ in $L^2(\mathbb{R}^n).$ Since $\mathcal{K}$ is non-zero, the continuous function $\varphi$ is nonzero, too. So, there  is such a constant $c>0$ that the set $\{s\in U: |\varphi(s)|>c\}$ contains some open ball $S.$ Then the restriction of $M_\varphi$ to the space $L^2(S)$ is an invertible and compact operator,
which is contrary to the fact that $L^2(S)$ is infinite dimensional.

Now let  $1\leq p< 2.$  It follows from (1) then
$M_{\psi}P_{L^q(S)}\mathcal{M}\mathcal{K}=P_{L^q(S)}\mathcal{M}$
where $\psi=(1/\varphi)|S$ and $P_{L^q(S)}f:=f|L^q(S)$ a projection. Let $T=P_{L^q(S)}\mathcal{M}.$ Passing to conjugates gives $\mathcal{K}^*T^*M_{\psi}^*=T^*.$ By the theorem 1 from \cite{FNR3} this implies that the operator $T^*=\mathcal{M}^*P_{L^q(S)}^*$ has finite rank. Note that the conjugate $P_{L^q(S)}^*$ is the operator of natural embedding $L^p(S)\subset L^p(\mathbb{R}^n).$ We clame also that for $g\in L^p(\mathbb{R}^n), x\in U$
$$
(\mathcal{M}^*g)(x)=\frac{1}{(2\pi)^{n/2}}\int_{\mathbb{R}^n}|x|^{-\frac{1}{q}+is}g(s)ds.
$$
First note that the right-hand side here is a bounded operator from $L^p(\mathbb{R}^n)$ to $L^q(U).$ As above this can be  obtained from the Hausdorff–Young inequality  by using an exponential change of variables. Then (and by lemma 1) both sides of the equality $(f\in L^p(U), g\in  L^q(\mathbb{R}^n))$
$$
\int_{\mathbb{R}^n}(\mathcal{M}f)(s)g(s)ds=\int_{U}f(x)(\mathcal{M}^*g)(x)dx
$$
 that we want to prove exist as Lebesgue integrals. So we may apply Fubini’s theorem
 as follows:
 $$
 \int_{U}f(x)(\mathcal{M}^*g)(x)dx=\int_{U}f(x)\frac{1}{(2\pi)^{n/2}}\int_{\mathbb{R}^n}|x|^{-\frac{1}{q}+is}g(s)dsdx=
 $$
$$
\int_{\mathbb{R}^n}\frac{1}{(2\pi)^{n/2}}\int_{U}f(x)|x|^{-\frac{1}{q}+is}dxg(s)ds=\int_{\mathbb{R}^n}(\mathcal{M}f)(s)g(s)ds.
$$
  It was shown above that the restriction of the operator $\mathcal{M}^*$  to $L^p(S)$ has finite rank. Since  $\mathcal{M}^*$ can be easily reduced to the Fourier transform, this contradicts to the Paley-Wiener Theorem about the image of $L^2(S)$ ($L^2(S)\subset L^p(S)$) under the Fourier transform (see, e. g., \cite[Theorem III.4.9]{SW}).

Finally, if  $2<p< \infty$ one cat use duality arguments. Indeed, by lemma 2
the adjoint operator $\mathcal{H}^*=\mathcal{H}_{\Phi, A'}^*$ (as an operator in $L^q(\mathbb{R}^n)$) is also of Hausdorff type. More precisely, it  equals to $\mathcal{H}_{\Psi, B},$ where $\Psi(u)=\Phi(u)|\det A(u)'^{-1}|=\Phi(u)/a(u),$ $B(u)= A(u)'^{-1}=\mathrm{diag}[1/a_1(u),\dots,1/a_n(u)].$ It is easy to verify that $\mathcal{H}_{\Psi, B}$ satisfies all the conditions of theorem 2 (with $q,$ $\Psi,$ and $B$ in place of $p,$ $\Phi,$ and $A$ respectively).
 Since  $1<q<2$, the   operator $\mathcal{H}_{\Psi, B}$ is  noncompact in $L^q(\mathbb{R}^n)$ and the result follows.

\

\textbf{Corollary 5.} \textit{Let $\forall j a_j(u)=b(u),$ where $b(u)$ is $\mu$-measurable and positive and $b(u)^{-n/p}\Phi(u)\in L^1(\Omega),$ $1\leq p<\infty.$ Then the corresponding Hausdorff operator $\mathcal{H}_\Phi$ is  noncompact in $L^p(\mathbb{R}^n)$ provided it is non-zero.
}

\

 Adolf R. Mirotin,

\textsc{Department of Mathematics and   Programming  Technologies,
F. Skorina Gomel State University}, 246019, Sovietskaya,
104, Gomel, Belarus

\end{document}